\documentclass[reqno,12pt]{amsart}
\usepackage{amssymb}
\usepackage{euscript}

\newcommand{\dfr}[2]{\dfrac{#1}{#2}}%ams
\newcommand{\cd}{\cdot}

\newcommand{\dsum}{\displaystyle \sum}

\newcommand{\abs}[1]{\lvert{#1}\rvert}

\newcommand{\qq}{\qquad}
\theoremstyle{plain}
\newtheorem{theorem}{Theorem}[section]

\newtheorem{lemma}[theorem]{Lemma}
\newtheorem{proposition}[theorem]{Proposition}
\theoremstyle{remark}
\newtheorem{remark}[theorem]{Remark}

\theoremstyle{definition}

\numberwithin{equation}{section}
\newcommand{\al}{\alpha}
\newcommand{\be}{\beta}
\newcommand{\bbe}{\boldsymbol{\beta}}

\newcommand{\om}{\omega}
\newcommand{\tom}{\tilde{\omega}}

\newcommand{\ep}{\epsilon}

\newcommand{\la}{\langle}

\newcommand{\ra}{\rangle}
\newcommand{\Vir}{\mathrm{Vir}}

\newcommand{\Z}{\mathbb{Z}}
\newcommand{\C}{\mathbb{C}}
\newcommand{\R}{\mathbb{R}}

\newcommand{\Q}{\mathbb{Q}}

\newcommand{\ba}{\mathbf{a}}

\DeclareMathOperator{\supp}{supp}
\DeclareMathOperator{\End}{End}

\begin{document}
\title[Vertex operator algebra and McKay's observation]
{Vertex operator algebras, extended $E_8$ diagram, and McKay's
observation on the Monster simple group}

\author[C.H. Lam]{Ching Hung Lam $^\dagger$}
\address[C.H.  Lam ]{Department of Mathematics, National Cheng Kung University,
Tainan, Taiwan 701}%
\email{chlam@mail.ncku.edu.tw}

\author[H. Yamada]{Hiromichi Yamada $^\ddagger$}
\address[H. Yamada]{
Department of Mathematics, Hitotsubashi University, Kunitachi,
Tokyo 186-8601,  Japan} \email{yamada@math.hit-u.ac.jp}

\author[H. Yamauchi]{  Hiroshi Yamauchi}
\address[H. Yamauchi]{Graduate School of Mathematics, University of
Tsukuba, Ibaraki 305-8571, Japan}

\email{hirocci@math.tsukuba.ac.jp}

\thanks{$^\dagger$ Partially supported by NSC grant 91-2115-M-006-014 of Taiwan,
R.O.C.\\
\hspace*{3mm} $^\ddagger$ Partially supported by JSPS Grant-in-Aid for
Scientific Research No. 15540015}%

\subjclass{17B68, 17B69, 20D08}

\begin{abstract}
We study McKay's observation on the Monster simple group, which
relates the $2A$-involutions of the Monster simple group to the
extended $E_8$ diagram, using the theory of vertex operator
algebras (VOAs). We first consider the sublattices $L$ of the
$E_8$ lattice obtained by removing one node from the extended
$E_8$ diagram at each time. We then construct a certain coset (or
commutant) subalgebra $U$ associated with $L$ in the lattice VOA
$V_{\sqrt{2}E_8}$. There are two natural conformal vectors of
central charge $1/2$ in $U$ such that their inner product is
exactly the value predicted by Conway \cite{c1}. The Griess
algebra of $U$ coincides with the algebra described in \cite[Table
3]{c1}. There is a canonical automorphism of $U$ of order
$|E_8/L|$. Such an automorphism can be extended to the Leech
lattice VOA $V_\Lambda$ and it is in fact a product of two
Miyamoto involutions. In the sequel \cite{lyy} to this article we
shall develop the representation theory of $U$. It is expected
that if $U$ is actually contained in the Moonshine VOA
$V^\natural$, the product of two Miyamoto involutions is in the
desired conjugacy class of the Monster simple group.
\end{abstract}
\maketitle

\section{Introduction}
The Moonshine vertex operator algebra $V^{\natural}$ constructed
by Frenkel-Lepowsky-Meurman \cite{flm} is one of the most
important examples of vertex operator algebras (VOAs). Its full
automorphism group is the Monster simple group. The weight $2$
subspace $V^{\natural}_2$ of $V^{\natural}$ has a structure of
commutative non-associative algebra which coincides with the
$196884$-dimensional algebra investigated by Griess \cite{griess}
in his construction of the Monster simple group (see also
Conway\cite{c1}). The structure of this algebra, which is called
the Monstrous Griess algebra, has been studied by group theorists.
It is well known \cite{c1} that each $2A$-involution $\phi$ of the
Monster simple group uniquely defines an idempotent $e_{\phi}$
called an axis in the Monstrous Griess algebra. Moreover, the
inner product $\la e_{\phi}, e_{\psi}\ra $ of any two axes
$e_{\phi}$ and $e_{\psi}$ is uniquely determined by the conjugacy
class of the product $\phi\psi$ of $2A$-involutions. Actually,
$2A$-involutions of the Monster simple group satisfy a
6-transposition property, that is, $|\,\phi\psi\,|\leq 6$ for any
two $2A$-involutions $\phi$ and $\psi$. In addition, the conjugacy
class of $\phi\psi$ is one of $1A$, $2A$, $3A$, $4A$, $5A$, $6A$,
$4B$, $2B$, or $3C$.

John McKay \cite{McK} observed that there is an interesting
correspondence with the extended $E_8$ diagram. Namely, one can
assign $1A, 2A, 3A, 4A, 5A, 6A,  4B, 2B$, and $3C$ to the nodes of
the extended $E_8$ diagram as follows (cf. Conway \cite{c1},
Glauberman and Norton \cite{GN}):

\begin{equation}\label{mckay}
\begin{array}{l}
  \hspace{184pt} 3C\quad   \frac{1}{2^8}\\
  \hspace{186.0pt}\circ \vspace{-6pt}\\
 \hspace{187.4pt}| \vspace{-6pt}\\
 \hspace{187.4pt}| \vspace{-6pt}\\
  \hspace{6pt} \circ\hspace{-5pt}-\hspace{-7pt}-\hspace{-5pt}-\hspace{-5pt}-
  \hspace{-5pt}-\hspace{-5pt}\circ\hspace{-5pt}-\hspace{-5pt}-
  \hspace{-5pt}-\hspace{-6pt}-\hspace{-7pt}-\hspace{-5pt}\circ
  \hspace{-5pt}-\hspace{-5.5pt}-\hspace{-5pt}-\hspace{-5pt}-
  \hspace{-7pt}-\hspace{-5pt}\circ\hspace{-5pt}-\hspace{-5.5pt}-
  \hspace{-5pt}-\hspace{-5pt}-\hspace{-7pt}-\hspace{-5pt}\circ
  \hspace{-5pt}-\hspace{-6pt}-\hspace{-5pt}-\hspace{-5pt}-
  \hspace{-7pt}-\hspace{-5pt}\circ\hspace{-5pt}-\hspace{-5pt}-
  \hspace{-6pt}-\hspace{-6pt}-\hspace{-7pt}-\hspace{-5pt}\circ
  \hspace{-5pt}-\hspace{-5pt}-
  \hspace{-6pt}-\hspace{-6pt}-\hspace{-7pt}-\hspace{-5pt}\circ
  \vspace{-6.2pt}\\
  \vspace{-6pt} \\
  1A\hspace{23pt} 2A\hspace{23  pt} 3A\hspace{22pt} 4A\hspace{21pt} 5A\hspace{21pt}
  6A\hspace{20pt} 4B\hspace{19pt} 2B\\
  \hspace{4pt} \frac{1}{4} \hspace{30pt} \frac{1}{32}\hspace{26  pt} \frac{13}{2^{10}}
  \hspace{22pt} \frac{1}{2^7} \hspace{25pt} \frac{3}{2^{9}}\hspace{23pt}
  \frac{5}{2^{10}}\hspace{23pt} \frac{1}{2^8}\hspace{23pt} 0 \\
  \
\end{array}
\end{equation}
where the numerical labels are equal to the multiplicities of the
corresponding simple roots in the highest root and the numbers
behind the labels denote the inner product $\la 2e_\phi,
2e_\psi\ra$ of $2e_\phi$ and $2e_\psi$.

%The Griess subalgebras generated by two axes $e_{\theta}$ and
%$e_{\phi}$ are determined by Conway\cite{c1}.

On the other hand, from the point of view of VOAs, Miyamoto
\cite{m1,m4} showed that an axis is essentially a half of a
conformal vector $e$ of central charge 1/2 which generates a
Virasoro VOA $\Vir(e)\cong L(1/2, 0)$ inside the Moonshine VOA
$V^{\natural}$. Moreover, an involutive automorphism $\tau_e$ can
be defined by
\[
\tau_e=
  \left\{
    \begin{array}{ll}
      1 & \mbox{on } W_0\oplus W_{1/2},\\
    -1 &  \mbox{on } W_{1/16},
    \end{array}
  \right.
\]
where $W_{h}$ denotes the sum of all irreducible $\Vir(e)$-modules
isomorphic to $L(1/2, h)$ inside $V^{\natural}$. In fact, $\tau_e$
is always of class $2A$ for any conformal vector $e$ of central
charge $1/2$ in $V^\natural$.

In this article, we try to give an interpretation of the McKay
diagram \eqref{mckay} using the theory of VOAs. We first observe
that there is a conformal vector $\hat{e}$ of central charge $1/2$
in the lattice VOA $V_{\sqrt{2}E_8}$ which is fixed by the action
of the Weyl group of type $E_8$.  Let $\Phi$ be the root system
corresponding to the Dynkin diagram obtained by removing one node
from the extended $E_8$ diagram and $L=L(\Phi)$ the root lattice
associated with $\Phi$. Then the Weyl group $W(\Phi)$ of $\Phi$
and the quotient group $E_8/L$ both act naturally on
$V_{\sqrt{2}E_8}$ and their actions commute with each other. The
action of the quotient group $E_8/L$ can be extended to the Leech
lattice VOA $V_\Lambda$.

The main idea is to construct certain vertex operator subalgebras
$U$ of the lattice VOA $V_{\sqrt{2}E_8}$ corresponding to the nine
nodes of the McKay diagram. In each case, $U$ is constructed as a
coset (or commutant) subalgebra of $V_{\sqrt{2}E_8}$ associated
with $\Phi$. In fact, $U$ is chosen so that the Weyl group
$W(\Phi)$ acts trivially on it. We show that in each of the nine
cases $U$ always contains $\hat{e}$ and another conformal vector
$\hat{f}$ of central charge $1/2$ such that the inner product $\la
\hat{e}, \hat{f}\ra$ is exactly the value listed in the McKay
diagram. Both of $\hat{e}$ and $\hat{f}$ are fixed by the Weyl
group $W(\Phi)$. Thus the Miyamoto involutions $\tau_{\hat{e}}$
and $\tau_{\hat{f}}$ commute with the action of $W(\Phi)$.
Furthermore, the quotient group $E_8/L$ naturally induces some
automorphism of $U$ of order $n = |E_8/L|$, which is identical
with the numerical label of the corresponding node in the McKay
diagram. Such an automorphism can be extended to the Leech lattice
VOA $V_\Lambda$ and it is in fact a product
$\tau_{\hat{e}}\tau_{\hat{f}}$ of two Miyamoto involutions
$\tau_{\hat{e}}$ and $\tau_{\hat{f}}$.

In the sequel \cite{lyy} to this article we shall study the
properties of the coset subalgebra $U$ in detail. Except the $4A$
case, $U$ always contains a set of mutually orthogonal conformal
vectors such that their sum is the Virasoro element of $U$ and the
central charge of those conformal vectors are all coming from the
unitary series
\begin{equation*}
c=c_m=1-\frac{6}{(m+2)(m+3)}, \qquad m=1,2,3, \dots.
\end{equation*}

Such a conformal vector generates a Virasoro VOA isomorphic to
$L(c_m,0)$ inside $U$. The structure of $U$ as a module for a
tensor product of those Virasoro VOA is determined.

In the $4A$ case, $U$ is isomorphic to the fixed point subalgebra
$V_\mathcal{N}^+$ of $\theta$ for some rank two lattice
$\mathcal{N}$, where $\theta$ is an automorphism of
$V_\mathcal{N}$ induced from the $-1$ isometry of the lattice
$\mathcal{N}$.

The VOA $U$ is generated by $\hat{e}$ and $\hat{f}$. As a
consequence we know that every element of $U$ is fixed by the Weyl
group $W(\Phi)$. The weight $1$ subspace $U_1$ of $U$ is $0$. The
Griess algebra $U_2$ of $U$ is also generated by $\hat{e}$ and
$\hat{f}$ and it has the same structure as the algebra studied in
Conway \cite[Table 3]{c1}. The automorphism group of $U$ is a
dihedral group of order $2n$ except the cases for $1A$, $2A$, and
$2B$. It is a trivial group in the $1A$ case, a symmetric group of
degree $3$ in the $2A$ case, and of order $2$ in the $2B$ case.
Furthermore, we shall discuss the rationality of $U$ and the
classification of irreducible modules. The product
$\tau_{\hat{e}}\tau_{\hat{f}}$ of two Miyamoto involutions should
be in the desired conjugacy class of the Monster simple group,
provided that the Moonshine VOA $V^\natural$ contains a subalgebra
isomorphic to $U$.

Further mysteries concerning the McKay diagram can be found in
Glauberman and Norton \cite{GN}. Among other things, some relation
between the Weyl group $W(\Phi)$ and the centralizer of a certain
subgroup generated by two $2A$-involutions and one $2B$-involution
in the Monster simple group was discussed. That every element of
$U$ is fixed by $W(\Phi)$ seems quite suggestive.

Let us recall some terminology (cf. \cite{flm}). A VOA is a
$\Z$-graded vector space $V = \oplus_{n \in \Z} V_n$ with a linear
map $Y(\,\cdot\,,z): V \rightarrow (\End V)[[z,z^{-1}]]$ and two
distinguished vectors; the vacuum vector $\mathbf{1} \in V_0$ and
the Virasoro element $\omega \in V_2$ which satisfy certain
conditions. For any $v \in V$, $Y(v,z) = \sum_{n \in \Z} v_n
z^{-n-1}$ is called a vertex operator and $v_n \in \End V$ a
component operator. Each homogeneous subspace $V_n$ is the
eigenspace for the operator $L(0) = \omega_1$ with eigenvalue $n$.
The eigenvalue for $L(0)$ is called a weight. Suppose $V =
\oplus_{n=0}^{\infty} V_n$ with $V_0 = \C\mathbf{1}$ and $V_1 =
0$. For $u, v \in V_2$, one can define a product $u \cdot v$ by
$u_1 v$ and an inner product $\la u,v\ra$ by $u_3 v = \la u,v\ra
\mathbf{1}$. The inner product is invariant, that is, $\la u_1 v,
w\ra = \la v, u_1 w\ra$ for $u,v,w \in V_2$ (cf. \cite[Section
8.9]{flm}). With the product and the inner product $V_2$ becomes
an algebra, which is called the Griess algebra of $V$.

The organization of the article is as follows. In Section 2 we
review some notation for lattice VOAs from \cite{flm} and certain
conformal vectors in the lattice VOA $V_{\sqrt{2}R}$ given by
\cite{dlmn}, where $R$ is a root lattice of type $A$, $D$, or $E$.
Moreover, we study some highest weight vectors in irreducible
modules of $V_{\sqrt{2}R}$ with respect to those conformal
vectors. In Section 3 we consider the sublattice $L$ of $E_8$ and
define the coset subalgebra $U$ and two conformal vectors
$\hat{e}$ and $\hat{f}$ of central charge $1/2$. We calculate the
inner product $\la\hat{e},\hat{f}\ra$ and verify that it is
identical with the value given in the McKay diagram. A canonical
automorphism $\sigma$ of order $n=|E_8/L|$ induced by the quotient
group $E_8/L$ is also discussed. Then in Section 4 we consider an
embedding of an orthogonal sum $\sqrt{2}{E_8}^3$ of three copies
of $\sqrt{2}{E_8}$ into the Leech lattice $\Lambda$ and show that
the product $\tau_{\hat{e}}\tau_{\hat{f}}$ of two Miyamoto
involutions $\tau_{\hat{e}}$ and $\tau_{\hat{f}}$ is of order $n$
as an automorphism of $V_{\Lambda}$. Finally, in Section 5 we give
an explicit correspondence between the Griess algebra $U_2$ of $U$
and the algebra in Conway \cite[Table 3]{c1}.

The authors thank Masaaki Kitazume and Masahiko Miyamoto for
stimulating discussions and Kazuhiro Yokoyama for helping them to
compute the conformal vectors for the cases of $5A$ and $6A$ by a
computer algebra system Risa/Asir.

\section{Conformal vectors in lattice VOAs} \label{S2}
In this section, we review the construction of certain conformal
vectors in the lattice VOA $V_{\sqrt{2}R}$ from \cite{dlmn}, where
$R$ is a root lattice of type $A_n$, $D_n$, or $E_n$. The notation
for lattice VOAs here is standard (cf. \cite{flm}). Let $N$ be a
positive definite even lattice with inner product $\la\,\cdot\, ,
\,\cdot\,\ra$. Then the VOA $V_N$ associated with $N$ is defined
to be $M(1) \otimes \C\{N\}$. More precisely, let $\mathfrak{h} =
\C \otimes_{\Z} N$ be an abelian Lie algebra and
$\hat{\mathfrak{h}} = \mathfrak{h} \otimes \C [t,t^{-1}] \oplus \C
K$ its affine Lie algebra. Then $M(1) = \C[\alpha(n)\,|\,\alpha
\in \mathfrak{h}, n < 0]\cdot 1$ is the unique irreducible
$\hat{\mathfrak{h}}$-module such that $\alpha(n) \cdot 1 = 0$ for
$\alpha \in \mathfrak{h}$, $n \ge 0$ and $K=1$, where $\alpha(n) =
\alpha \otimes t^n$. Moreover, $\C\{N\}$ denotes a twisted group
algebra of the additive group $N$. In the case for $N =
\sqrt{2}R$, the twisted group algebra $\C\{\sqrt{2}R\}$ is
isomorphic to the ordinary group algebra $\C[\sqrt{2}R]$ since
$\sqrt{2}R$ is a doubly even lattice. The standard basis of
$\C[\sqrt{2}R]$ is denoted by $e^{\sqrt{2}\al}$, $\al \in R$. Then
the vacuum vector $\mathbf{1}$ is $1 \otimes e^0$.

Let $\Phi$ be the root system of $R$ and $\Phi^+$ and $\Phi^-$ the
set of all positive roots and negative roots, respectively. Then
$\Phi=\Phi^+\cup \Phi^-= \Phi^+\cup ( -\Phi^+)$. The Virasoro
element $\omega$ of $V_{\sqrt{2}R}$ is given by
\[
\omega =\om (\Phi)=\frac{1}{2h}\sum_{\al\in \Phi^+} \al(-1)^2\cdot
1,
\]
where $h$ is the Coxeter number of $\Phi$. Now define
\begin{equation}\label{cv1}
\begin{split}
s&=s(\Phi)=\frac{1}{2(h+2)}\sum_{\al\in \Phi^+}\left(
\al(-1)^2\cdot
1 -2(e^{\sqrt{2}\al}+ e^{-\sqrt{2}\al})\right),\\
\tilde{\om}&=\tilde{\om}(\Phi)= \om - s.
\end{split}
\end{equation}
It is shown in \cite{dlmn} that $\tilde{\om}$ and $s$ are mutually
orthogonal conformal vectors, that is, $\tilde{\om}_1\tilde{\om}=
2\tilde{\om}, s_1 s=2s$, and $\tilde{\om}_1s=0$. The central
charge of $\tilde{\om}$ is $2n/(n+3)$ if $R$ is of type $ A_n$,
$1$ if $R$ is of type $ D_n$ and $6/7, 7/10$ and $1/2$ if $R$ is
of type $E_6, E_7$ and $E_8$, respectively.

Let $W(\Phi)$ be the Weyl group of $\Phi$. Any element $g\in
W(\Phi)$ induces an automorphism of the lattice $R$ and hence it
defines an automorphism of the  VOA $V_{\sqrt{2}R}$ by
\[
g(u\otimes e^{\sqrt{2}\al})= gu\otimes e^{\sqrt{2} g\al} \qquad
\text{ for } \quad u\otimes e^{\sqrt{2}\al}\in M(1)\otimes
e^{\sqrt{2}\al}\subset V_{\sqrt{2}R}.
\]
Note that both $s$ and $\tilde{\om}$ are fixed by the Weyl group
$W(\Phi)$.
%In fact, $\om,\tom$ and $s$ are the only conformal
%vectors in $V_{\sqrt{2}R}$ which are fixed by the Weyl group.

We shall study certain highest weight vectors with respect to the
subalgebra $\Vir(s)\otimes \Vir(\tilde{\om})$, where $\Vir(s)$ and
$\Vir(\tilde{\om})$ denote the Virasoro VOAs generated by the
conformal vectors $s$ and $\tilde{\om}$, respectively.

Let $R^* = \{ \al \in \Q \otimes_{\Z} R\,|\,\la\al,R\ra \subset \Z\}$
be the dual lattice of $R$.

\begin{lemma}
Let $R$ be a root lattice of type $A$, $D$, or $E$ and $\gamma+R$
a coset of $R$ in $R^*$. Let $k=\min\{ \la \al,\al\ra | \al\in
\gamma+ R\}$. For any $\eta\in \gamma+R$ with $\la \eta,\eta\ra
=k$, we define
\[
X_\eta=\{(\al,\be)\in R\times (\gamma+R)\ |\ \la \al,\al\ra =2,
\la\be,\be \ra=k \text{ and } \al+\be=\eta\}.
\]
Then $|X_\eta|= kh$, where $h$ is the Coxeter number of $R$.
\end{lemma}

\medskip

\begin{proof}
The proof is just by direct verification.  We only discuss the
case for $R=A_n$. The other cases can be proved similarly.

Let $R=A_n$. Then the Coxeter number $h$ is $n+1$ and the roots of
$A_n$ are  given by the vectors in the form $\pm ( 1,-1, 0^{n-1})
\in \R^{n+1}$, that is, the vectors whose one entry is $\pm 1$,
another entry is $\mp 1$, and the remaining $n-1$ entries are $0$.
Let $\mu= \frac{1}{n+1} (1,\dots,1, -n)$. Then $\mu+ R$ is a
generator of the group $R^*/R$. Denote $\gamma= j\mu$ for $j=0,
\dots, n$. Then
\[
k=\min \{\la\al,\al\ra | \al\in \gamma+R\} = \frac{j(n+1-j)}{n+1},
\]
and the elements of square norm $k$ in $\gamma+R$ are of the form
\[
\frac{1}{n+1} ( j^{n+1-j}, (-n-1+j)^j).
\]

Now it is easy to see that $|X_{\eta}|= (n+1-j)j=kh$ for any $\eta
$ with $ \la \eta,\eta \ra =k$.
\end{proof}

\begin{proposition}\label{hw}
Let $\gamma+R$ be a coset of $R$ in $R^*$ and $k=\min\{ \la
\al,\al\ra | \al\in \gamma+ R\}$. Define
\[
v=\sum_{{\al\in \gamma+R}\atop \la \al,\al\ra=k} e^{\sqrt{2}\al}
\in V_{\sqrt{2}(\gamma+R)}.
\]
Then $v$ is a highest weight vector of highest weight $(0,k)$ in
$V_{\sqrt{2}(\gamma+ R)}$ with respect to  $\Vir(s)\otimes
\Vir(\tilde{\om})$, that is, $s_j v = \tilde{\om}_j v=0$ for all
$j \ge 2$, $s_1 v=0$, and $\tilde{\om}_1 v=kv$.
\end{proposition}

\begin{proof}
Since $k$ is the minimum weight of $V_{\sqrt{2}(\gamma+ R)}$, it
is clear that $s_j v= \tom_j v =0$ for all $j\geq 2$. Since
$\omega_1 v = kv$, it suffices to show that $s_1 v=0$. By the
definition \eqref{cv1} of $s$ and the above lemma, we have
\begin{align*}
s_1 v& =\frac{1}{2(h+2)}\sum_{\al\in \Phi^+}\left( \al(-1)^2\cdot
1 -2(e^{\sqrt{2}\al}+ e^{-\sqrt{2}\al})\right)_1 v\\
& = \Big( \frac{h}{h+2} \omega - \frac{1}{h+2}\sum_{\al\in
\Phi^+}(e^{\sqrt{2}\al}+
e^{-\sqrt{2}\al})\Big)_1 v\\
& =  \frac{hk}{h+2}\, v- \frac{hk}{h+2}\, v =0.
\end{align*}
Hence the assertion holds.
\end{proof}

\section{ Extended $E_8$ diagram and sublattices of the root
lattice $E_8$}\label{S3}

In this section,  we consider  certain sublattices of the root
lattice $E_8$ by using the extended $E_8$ diagram
\begin{equation}\label{ext}
\begin{array}{l}
  \hspace{183pt} \al_8\\
  \hspace{186.0pt}\circ \vspace{-6pt}\\
 \hspace{187.4pt}| \vspace{-6pt}\\
 \hspace{187.4pt}| \vspace{-6pt}\\
  \hspace{6pt} \circ\hspace{-5pt}-\hspace{-7pt}-\hspace{-5pt}-\hspace{-5pt}-
  \hspace{-5pt}-\hspace{-5pt}\circ\hspace{-5pt}-\hspace{-5pt}-
  \hspace{-5pt}-\hspace{-6pt}-\hspace{-7pt}-\hspace{-5pt}\circ
  \hspace{-5pt}-\hspace{-5.5pt}-\hspace{-5pt}-\hspace{-5pt}-
  \hspace{-7pt}-\hspace{-5pt}\circ\hspace{-5pt}-\hspace{-5.5pt}-
  \hspace{-5pt}-\hspace{-5pt}-\hspace{-7pt}-\hspace{-5pt}\circ
  \hspace{-5pt}-\hspace{-6pt}-\hspace{-5pt}-\hspace{-5pt}-
  \hspace{-7pt}-\hspace{-5pt}\circ\hspace{-5pt}-\hspace{-5pt}-
  \hspace{-6pt}-\hspace{-6pt}-\hspace{-7pt}-\hspace{-5pt}\circ
  \hspace{-5pt}-\hspace{-5pt}-
  \hspace{-6pt}-\hspace{-6pt}-\hspace{-7pt}-\hspace{-5pt}\circ
  \vspace{-6.2pt}\\
  \vspace{-6pt} \\
  \hspace{6pt}\al_0\qq \alpha_1\qq \alpha_2\qq \alpha_3\qq
  \alpha_4\qq \alpha_5\qq \alpha_6 \qq \alpha_7
\end{array}
\end{equation}
where $\al_1,\al_2,\dots,\al_8$ are the simple roots of $E_8$ and
\begin{equation}\label{eq:A_0}
\al_0+2\al_1+3\al_2+4\al_3+ 5\al_4+6\al_5+4\al_6+2\al_7+3\al_8 =
0.
\end{equation}
Thus $\la \al_i,\al_i\ra = 2$, $0 \le i \le 8$. Moreover, for
$i\ne j$, $\la \al_i, \al_j \ra = -1$ if the nodes $\al_i$ and
$\al_j$ are connected by an edge and $\la \al_i, \al_j \ra = 0$
otherwise. Note that $-\alpha_0$ is the highest root.

For any $i=0,1,\dots, 8$, let $L(i)$ be the sublattice generated
by $\al_j, 0 \le j \le 8, j \ne i$. Then $L(i)$ is a rank $8$
sublattice of $E_8$. In fact, $L(i)$ is the lattice associated
with the Dynkin diagram obtained by removing the corresponding
node $\al_i$ from the extended $E_8$ diagram \eqref{ext}. Note
that the index $| E_8/L(i) |$ is equal to $n_i$, where $n_i$ is
the coefficient of $\al_i$ in the left hand side of
\eqref{eq:A_0}. Actually, we have
\begin{gather}\label{eq:DEC}
\begin{aligned}
    L(0)&\cong E_8,\\
    L(3)&\cong A_3\oplus D_5,\\
    L(6)&\cong A_7\oplus A_1,
\end{aligned}
\qquad
\begin{aligned}
    L(1)&\cong A_1\oplus E_7,\\
    L(4)&\cong A_4\oplus A_4,\\
    L(7)&\cong D_8,
\end{aligned}
\qquad
\begin{aligned}
    L(2)&\cong A_2\oplus E_6, \\
    L(5)&\cong A_5\oplus A_2\oplus A_1,\\
    L(8)&\cong A_8.
\end{aligned}
\end{gather}

\begin{remark}
If $n_i$ is not a prime, there is an intermediate sublattice as
follows.
\begin{align*}
& A_3 \oplus D_5 \subset D_8 \subset E_8,\\
& A_5 \oplus A_2 \oplus A_1 \subset A_2 \oplus E_6 \subset E_8,
\qquad A_5 \oplus A_2 \oplus A_1 \subset A_1 \oplus E_7 \subset
E_8,\\
& A_7 \oplus A_1 \subset A_1 \oplus E_7 \subset E_8.
\end{align*}

There are corresponding power maps between conjugacy classes of
the Monster simple group, namely,
\begin{equation*}
(4A)^2 = 2B, \qquad (6A)^2 = 3A, \qquad (6A)^3 = 2A, \qquad (4B)^2
= 2A,
\end{equation*}
where $(mX)^k = nY$ means that the $k$-th power $g^k$ of an
element $g$ in the conjugacy class $mX$ is in the conjugacy class
$nY$ (cf. \cite{atlas}).
\end{remark}

\subsection{Coset subalgebras of the lattice VOA
$V_{\sqrt{2}E_8}$}\label{coset}

We shall construct some VOAs $U$ corresponding to the nine nodes
of the McKay diagram (\ref{mckay}). In each case, we show that the
VOA $U$ contains some conformal vectors of central charge $1/2$
and the inner products among these conformal vectors are the same
as the numbers given in the McKay diagram.

Let us explain the details of our construction. First, we fix
$i\in \{0,1,\dots,8\}$ and denote $L(i)$ by $L$. In each case,
$|E_8/L|=n_i$ and $\al_i+L$ is a generator of the quotient group
$E_8/L$. Hence we have
\begin{equation}\label{eq:E8L}
E_8=L\cup (\al_i+L)\cup (2\al_i+L)\cup \cdots \cup (
(n_i-1)\al_i+L).
\end{equation}
Then the lattice VOA $V_{\sqrt{2}E_8}$ can be decomposed as
\[
V_{\sqrt{2}E_8}=V_{\sqrt{2}L}\oplus
V_{\sqrt{2}\alpha_i+\sqrt{2}L}\oplus \cdots\oplus
V_{\sqrt{2}(n_i-1)\alpha_i+\sqrt{2}L},
\]
where $V_{\sqrt{2}j\alpha_i+\sqrt{2}L}$, $j=0,1,\dots,n_i-1$, are
irreducible modules of $V_{\sqrt{2}L}$ (cf. \cite{d1}).

The quotient group $E_8/L$ induces an automorphism $\sigma$ of
$V_{\sqrt{2}E_8}$ such that
\begin{equation}\label{au}
\sigma(u)=\xi^j u\qquad \text{ for any } \quad u\in
V_{\sqrt{2}j\alpha_i+\sqrt{2}L},
\end{equation}
where $\xi=e^{2\pi \sqrt{-1}/ n_{i}}$ is a primitive $n_i$-th root
of unity. More precisely, let
\begin{equation}\label{adef}
\ba = \begin{cases}
\al_1 & \mbox{if } i=0,\\
-\frac{1}{i+1}(\al_0 + 2\al_1 + \cdots + i\al_{i-1}) &
\mbox{if } 1\le i\le 5,\\
-\frac{1}{8}(\al_0+2\al_1 + \cdots + 6\al_5+7\al_8) &
\mbox{if } i=6,\\
\frac{1}{2}(\al_6+\al_8) &
\mbox{if } i=7,\\
-\frac{1}{9}(\al_0+2\al_1 + \cdots + 8\al_7) & \mbox{if } i=8.
\end{cases}
\end{equation}
Then $\la\ba,\al_j\ra \in \Z$ for $0 \le j\le 8$ with $j\ne i$ and
$\la\ba,\al_i\ra \equiv -1/n_i \pmod \Z$. The automorphism
$\sigma: V_{\sqrt{2}E_8} \to V_{\sqrt{2}E_8}$ is in fact defined
by
\begin{equation}\label{sigmadef}
\sigma = e^{-\pi \sqrt{-1} \bbe(0)} \qquad \text{ with } \quad
\bbe=\sqrt{2} \ba.
\end{equation}
For $u\in M(1)\otimes e^{\al} \subset V_{\sqrt{2}E_8}$, we have
$\sigma(u)=e^{-\pi \sqrt{-1} \la\bbe, \al\ra} u$. Note that $\ba
+R$ is a generator of the quotient group $R^*/R$ for the cases
$i\ne 0,7$, where $R$ is an indecomposable component of the
lattice $L$ of type $A$ and $R^*$ is the dual lattice of $R$.

For any lattice VOA $V_N$ associated with a positive definite even
lattice $N$, there is a natural involution $\theta$ induced by the
isometry $\al \to -\al$ for $\al \in N$. If $N = \sqrt{2}E_8$,
which is doubly even, we may define $\theta: V_{\sqrt{2}E_8} \to
V_{\sqrt{2}E_8}$ by
\begin{equation}\label{eq:theta}
\alpha(-n)\rightarrow -\alpha(-n)\qquad \mbox{ and }\qquad
e^{\alpha}\rightarrow e^{-\alpha}
\end{equation}
for $\al \in \sqrt{2}E_8$ (cf. \cite{flm}). Then $\theta \sigma
\theta= \sigma^{-1}$ and the group generated by $\theta$ and
$\sigma$ is a dihedral group of order $2n_i$.

Let $R_1,\dots, R_l$ be the indecomposable components of the
lattice $L$ and $\Phi_1,\dots, \Phi_l$ the corresponding root
systems of $R_1,\dots, R_l$ (cf. \eqref{eq:DEC}). Then
$L=R_1\oplus \cdots\oplus R_l$ and
\begin{equation*}
V_{\sqrt{2}L}\cong V_{\sqrt{2}R_1}\otimes \cdots \otimes
V_{\sqrt{2}R_l},
\end{equation*}
(see \cite{fhl} for tensor products of VOAs). By \eqref{cv1}, one
obtains $2l$ mutually orthogonal conformal vectors
\begin{equation}\label{eq:CV}
s^k=s(\Phi_k),\quad \tilde{\om}^k=\tilde{\om}(\Phi_k),\quad k=1,
\ldots, l
\end{equation}
such that the Virasoro element $\om$ of $V_{\sqrt{2}L}$,
which is also the Virasoro element of $V_{\sqrt{2}E_8}$, can be
written as a sum of these conformal vectors
\begin{equation*}
\om=s^1+\cdots+ s^l +\tilde{\om}^1+\cdots+ \tilde{\om}^l.
\end{equation*}

Now we define $U$ to be a coset (or commutant) subalgebra
\begin{equation}\label{eq:UDEF}
U=\{v\in V_{\sqrt{2}E_8}\,| \, (s^k)_1 v=0 \text{ for all } k=1,
\dots, l\}.
\end{equation}
Note that $U$ is a VOA with the Virasoro element
$\om'=\tilde{\om}^1+\cdots+ \tilde{\om}^l$ and the automorphism
$\sigma$ defined by (\ref{au}) induces an automorphism of order
$n_i$ on $U$. By abuse of notation, we denote it by $\sigma$ also.

\begin{remark}
In \cite{ly3}, it is shown that $\{ v\in V_{\sqrt{2}A_n}\,|\,
s(A_n)_1 v=0\}$ is isomorphic to a parafermion algebra
$W_{n+1}(2n/(n+3))$ of central charge $2n/(n+3)$. Thus, if $L$ has
some indecomposable component of type $A_n$, then $U$ contains
some subalgebra isomorphic to a parafermion algebra. It is well
known \cite{zf1} that the parafermion algebra $W_{n+1}(2n/(n+3))$
possesses a certain $\Z_{n+1}$ symmetry in the fusion rules among
its irreducible modules. The automorphism $\sigma$ is in fact
related to such a symmetry. More details about the relation
between coset subalgebra $U$ and the parafermion algebra
$W_{n+1}(2n/(n+3))$ can be found in \cite{lyy}.
\end{remark}

\subsection{ Conformal vectors of central charge $1/2$}
Next, we shall study some conformal vectors in $V_{\sqrt{2}E_8}$.
We shall also show that the coset subalgebra $U$ always contains
some conformal vectors of central charge $1/2$. Moreover, the
inner products among these conformal vectors will be discussed.

Recall that the lattice $\sqrt{2}E_8$ can be constructed by using
the $[8,4,4]$ Hamming code $H_8$ and the Construction A (cf.
\cite{cs}). That means
\begin{equation}\label{eq:2E8}
\sqrt{2}E_8 = \left \{ (a_1, \dots, a_8)\in \Z^8\, |\ (a_1, \dots,
a_8)\in H_8 \mod 2\right \}.
\end{equation}

We denote the vectors $(0,0,0,0,0,0,0,0)$ and $ (1,1,1,1,1,1,1,1)$
by $\mathbf{0}$ and $\mathbf{1}$, respectively. For any $\gamma\in
H_8$, we define
\begin{align*}
X^0_\gamma &=\sum_{ {\al\equiv \gamma\,\mathrm{mod}\, 2} \atop
{\la \al,\al\ra=4}} (-1)^{\la \al, \mathbf{0} \ra/2 } e^\al=\sum_{
{\al\equiv \gamma\,\mathrm{mod}\, 2} \atop {\la
\al,\al\ra=4}} e^\al,\\
X^1_\gamma &=\sum_{ {\al\equiv \gamma\,\mathrm{mod}\, 2} \atop
{\la \al,\al\ra=4}} (-1)^{\la \al, \mathbf{1} \ra/2 } e^\al,
\end{align*}
and for any binary word $\delta\in {\Z_2}^8$, we define
\begin{equation*}
\hat{e}^{\epsilon}_\delta = \frac{1}{16} \om + \frac{1}{32}
\sum_{\gamma \in H_8} (-1)^{\la\delta, \gamma\ra}
X^{\epsilon}_\gamma, \qquad \epsilon = 0,1,
\end{equation*}
where $\om$ is the Virasoro element of the VOA $V_{\sqrt{2}E_8}$.
Note that $X^{\ep}_{\mathbf{1}}=0$ for any $\ep=0,1$ and that
$\hat{e}^\ep_\delta= \hat{e}^\ep_{\eta}$ if and only if $\eta \in
\delta+H_8$
\medskip

\begin{lemma}\label{hate}
For any $\ep=0,1$ and $\delta\in {\Z_2}^8$, $\hat{e}^\ep_\delta$
is a conformal vector of central charge $1/2$. The inner product
among them are as follows.
\[
\la \hat{e}^\ep_\delta,  \hat{e}^\ep _\eta\ra  =
\begin{cases}
\displaystyle
0 & \text{if } \delta+\eta \text{ is even}\\
1/32 &\text{if } \delta+\eta \text{ is odd }
\end{cases}
\]
for any $\eta \notin \delta+H_8$, and
\[
\la \hat{e}^0_\delta,  \hat{e}^1_\eta\ra  = 0
\]
for any $\delta,\eta \in {\Z_2}^8$.
\end{lemma}

\begin{proof}
We have
\begin{align*}
(X_\gamma^\epsilon)_1(X_\zeta^\epsilon)&= 4 X^\epsilon_{\gamma
+\zeta}\qquad \text{ if }\qquad \abs{\gamma +\zeta}=4,\\
(X_{\mathbf{0}}^\epsilon)_1(X_{\mathbf{0}}^\epsilon) &=
\dsum_{\alpha \equiv \mathbf{0}\, \mathrm{mod}\, 2 \atop \la
\alpha,\alpha\ra =4} \dfr{1}{2} \alpha (-1)^2 \cd 1.
\end{align*}
Moreover, for any $\gamma \in H_8$ with $\abs{\gamma}=4$,
\begin{equation*}
  (X^\epsilon_\gamma)_1(X^\epsilon_\gamma)
  + (X^\epsilon_{\mathbf{1}+\gamma})_1 (X_{\mathbf{1}+\gamma}^\epsilon)
  = \dsum_{\alpha \equiv \gamma \, \mathrm{mod}\, 2
    \atop \la \alpha,\alpha \ra =4} \dfr{1}{2} \alpha (-1)^2\cd 1
  + \dsum_{\alpha \equiv \mathbf{1} +\gamma \, \mathrm{mod}\, 2 \atop
    \la \alpha,\alpha \ra =4} \dfr{1}{2} \alpha (-1)^2\cd 1
  + 8 X^\epsilon_{\mathbf{0}} .
\end{equation*}
Note also that
\begin{equation*}
\sum_{\gamma \in H_8} \sum_{\alpha \equiv \gamma \, \mathrm{mod}\,
2 \atop \la \alpha,\alpha \ra=4} \dfr{1}{2} \alpha (-1)^2\cd 1
=\sum_{\beta \in \Phi (E_8)} \beta (-1)^2 \cd 1 =2\sum_{\beta\in
\Phi^+ (E_8)} \beta(-1)^2\cd 1.
\end{equation*}

In addition, we have
\begin{equation*}
\la X^\ep_\gamma, X^\ep_\zeta \ra =
\begin{cases}
16 & \text{if } \gamma= \zeta \text{ and } \la \gamma,\gamma
\ra\neq 8, \\
0 & \text{otherwise},
\end{cases}
\end{equation*}
\begin{equation*}\label{01}
\la X^0_\gamma, X^1_\zeta \ra =
\begin{cases}
-16 & \text{if } \gamma= \zeta =0, \\
0 & \text{otherwise}.
\end{cases}
\end{equation*}

Then since $\omega_1 \omega = 2\omega$ and $\la \omega,\omega\ra =
4$, it follows that
\begin{equation*}
\begin{split}
(\hat{e}^\ep_{\delta})_1\hat{e}^\ep_\delta &=\Big( \frac{1}{16}
\om +\frac{1}{32} \sum_{\gamma\in H_8} (-1)^{\la\delta, \gamma\ra}
X^\ep_\gamma\Big)_1
 \Big(\frac{1}{16} \om +\frac{1}{32} \sum_{\gamma\in
H_8} (-1)^{\la\delta, \gamma\ra} X^\ep_\gamma\Big)\\
&= \frac{1}{2^8} \times
2\om+2\times\frac{1}{16}\times\frac{1}{32}\times2\sum_{\gamma\in
H_8}
(-1)^{\la\delta, \gamma\ra} X^\ep_\gamma\\
& \qquad +
\frac{1}{2^{10}}\Big(\sum_{\beta\in\Phi^{+}(E_8)}2\beta(-1)^2\cdot
1 + 56\sum_{\gamma\in
H_8} (-1)^{\la\delta, \gamma\ra} X^\ep_\gamma\Big)\\
&= \frac{1}{8}\om+\frac{1}{16}\sum_{\gamma\in H_8}
(-1)^{\la\delta,\gamma\ra} X^\ep_\gamma= 2\hat{e}^\ep_\delta,
\end{split}
\end{equation*}
and
\[
\la\hat{e}_\delta^\ep, \hat{e}_\delta^\ep \ra =
\frac{1}{2^{8}}\times 4+\frac{1}{2^{10}}\times 240 = \frac{1}{4}.
\]
Hence $\hat{e}_\delta^\ep$ is a conformal vectors of central
charge $1/2$.

For any $\eta \notin \delta+H_8$, we calculate that
\begin{equation*}
\begin{split}
\la\hat{e}^\ep_\delta , \hat{e}^\ep_\eta\ra &=
\frac{1}{2^{8}}\times 4+ \frac{1}{2^{10}} \sum_{ \gamma\in H_8
} (-1)^{\la\delta+\eta, \gamma\ra} \la X^\ep_\gamma, X^\ep_\gamma\ra\\
&=
\begin{cases}
\frac{1}{64}+ \frac{1}{2^{10}}\times 16 \times (7-8) = 0 &
\text{ if } \delta+\eta \text { is even}, \\
\frac{1}{64}+ \frac{1}{2^{10}}\times 16 \times (8-7) =
\frac{1}{32} & \text{ if } \delta+\eta \text { is odd}.
\end{cases}
\end{split}
\end{equation*}
Note that there are exactly eight elements in $H_8$ which are
orthogonal to $\delta+\eta$. Note also that $\delta+\eta$ is
orthogonal to $(1,1,1,1,1,1,1,1)$  if and only if $\delta+\eta$ is
even.

Finally, for any $\delta, \eta \in {\Z_2}^8$ we obtain
\begin{equation*}
\la\hat{e}^0_\delta, \hat{e}^1_\eta\ra = \frac{1}{2^{8}}\times 4-
\frac{1}{2^{10}}\times 16=0.
\end{equation*}
\end{proof}

In Miyamoto \cite{m2}, certain conformal vectors of central charge
$1/2$ are constructed inside the Hamming code VOA. Our
construction of $\hat{e}^\ep_\delta$ is essentially a lattice
analogue of Miyamoto's construction. In fact, take
$\lambda_j=(0,\dots,2,\dots,0)\in \Z^8$ to be the element in
$\sqrt{2}E_8$ such that the $j$-th entry is 2 and all other
entries are zero. Then we have a set of $16$ mutually orthogonal
conformal vectors of central charge $1/2$ given by
\[
\om^\pm_{\lambda_j}=\frac{1}{16} \lambda_j(-1)^2\cdot 1 \pm
\frac{1}{4} (e^{\lambda_j}+e^{-\lambda_j}), \quad j = 1,2,\dots,8.
\]

A set of mutually orthogonal conformal vectors of central charge
$1/2$ whose sum is equal to the Virasoro element in a VOA is
called a Virasoro frame. Thus, $\{\om^\pm_{\lambda_j}\,|\,1 \le j
\le 8\}$ is a Virasoro frame of $V_{\sqrt{2}E_8}$. With respect to
this Virasoro frame, the lattice VOA $V_{\sqrt{2}E_8}$ is a code
VOA (cf. \cite{m2}). Let $V_{\sqrt{2}{E_8}}^+$ be the fixed point
subalgebra of $V_{\sqrt{2}E_8}$ under the automorphism $\theta$
(cf. \eqref{eq:theta}). Then $\om^\pm_{\lambda_j}\in
V_{\sqrt{2}{E_8}}^+$ and $V^+_{\sqrt{2}E_8}$ is isomorphic to a
code VOA $M_D$, where $D$ is the second order Reed-M\"{u}ller code
$RM(4,2)$ of length $16$. Note that $\dim RM(4,2)=11$ and the dual
code of $RM(4,2)$ is the first order Reed-M\"{u}ller code
$RM(4,1)$ with the generating matrix
\[
\left(%
\begin{array}{cccccccccccccccc}
 1&1&1&1&1&1&1&1&1&1&1&1&1&1&1&1\\
 1&1&1&1&1&1&1&1&0&0&0&0&0&0&0&0\\
 1&1&1&1&0&0&0&0&1&1&1&1&0&0&0&0\\
 1&1&0&0&1&1&0&0&1&1&0&0&1&1&0&0\\
 1&0&1&0&1&0&1&0&1&0&1&0&1&0&1&0
\end{array}%
\right)
\]

Let $H^+$ and $H^-$ be the subcodes of $D$ whose supports are
contained in the positions corresponding to
$\{\om^+_{\lambda_j}\,|\, 1 \le j \le 8\}$ and
$\{\om^-_{\lambda_j}\,|\, 1 \le j \le 8\}$, respectively. Then
$H^+$ and $H^-$ are both isomorphic to the $[8,4,4]$ Hamming code
$H_8$. The conformal vectors $\hat{e}^0_\delta$ and
$\hat{e}^1_\delta$ are actually the conformal vectors $s_\delta$
constructed by Miyamoto\,\cite{m2} using the Hamming code VOAs
$M_{H^+}$ and $M_{H^-}$, respectively.

\begin{proposition}\label{frame}
The set $\{\hat{e}^0_\delta, \hat{e}^1_\zeta\,|\, \delta, \zeta
\in {\Z_2}^8/H_8, \delta, \zeta \text {are even}\}$  is a Virasoro
frame of $V_{\sqrt{2}{E_8}}^+$. Moreover, $V^+_{\sqrt{2}E_8}\cong
M_{RM(4,2)}$ with respect to this frame, where $M_{RM(4,2)}$
denotes the code VOA associated with the second order
Reed-M\"{u}ller code $RM(4,2)$.
\end{proposition}

\begin{proof}
The first assertion follows from Lemma \ref{hate}. As mentioned
above, we know that $V_{\sqrt{2}E_8}^+\cong M_D$ with respect to
the frame $\{\om^\pm_{\lambda_j}\,|\, 1 \le j \le 8 \}$, where $D
\cong RM(4,2)$. It contains a subalgebra isomorphic to
$M_{H^+}\otimes M_{H^-}$. For convenience, we arrange the
positions of $\{\om^\pm_{\lambda_j}\}$ so that the support $\supp
H^+$ of $H^+$ is $(1^8,0^8)$ and the support $\supp H^-$ of $H^-$
is $(0^8,1^8)$. Let $\{\be_0, \be_1, \dots, \be_7\}$ with
$\be_0=0$ be a complete set of coset representatives of
$D/(H^+\oplus H^-)$. Then
\[
V_{\sqrt{2}E_8}^+ \cong M_{H^+\oplus H^-} \oplus \bigoplus_{i=1}^7
M_{\be_i+(H^+\oplus H^-)}.
\]

By a result of Miyamoto \cite{m2}, $M_{H^+\oplus H^-}$ is still
isomorphic to the code VOA $M_{H^+\oplus H^-}$ associated with
$H^+\oplus H^-$ with respect to the frame $\{\hat{e}^0_\delta,
\hat{e}^1_\zeta\,|\, \delta, \zeta \in {\Z_2}^8/H_8, \delta, \zeta
\text {are even}\}$. Moreover, we know that $(1^8,0^8)$ and
$(0^8,1^8)$ are contained in the dual code of $D$. Thus
$\la(1^8,0^8), \be_i\ra= \la(0^8,1^8), \be_i\ra=0$ for all $i$.
Let $\be^+$ and $\be^-$ be such that $\mathrm{supp} \be^+\subset
\mathrm{supp} H^+$, $\mathrm{supp} \be^-\subset \mathrm{supp}
H^-$, and $\be_i= \be^++\be^-$. Then $M_{\be_i+(H^+\oplus
H^-)}\cong M_{\be^++H^+}\otimes M_{\be^-+H^-}$ and both of
$M_{\be^++H^+}$ and $ M_{\be^-+H^-}$  are of integral weight.
Hence, by \cite{m2}, $M_{\be_i+(H^+\oplus H^-)}$ is again
isomorphic to $ M_{\be_i+(H^+\oplus H^-)}$ with respect to the
frame $\{\hat{e}^0_\delta, \hat{e}^1_\zeta\,|\, \delta, \zeta \in
{\Z_2}^8/H_8, \delta, \zeta \text {are even}\}$ and thus we still
have $ V^+_{\sqrt{2}E_8}\cong M_{D}$.
\end{proof}

Now let
\begin{equation}\label{eq:EF}
\begin{aligned}
\hat{e} &= \hat{e}^0_{\mathbf{0}}= \frac{1}{16}\om
+\frac{1}{32}\sum_{\al\in\Phi^+(E_8)}
(e^{\sqrt{2}\al}+e^{-\sqrt{2}\al}),\\
\hat{f} &= \sigma \hat{e},
\end{aligned}
\end{equation}
where $\sigma$ is the automorphism defined by \eqref{au}. These
conformal vectors of central charge $1/2$ play an important role
for the rest of the paper.

Let $\Phi$ be the root system of $L=L(i)$. Let $H_j = \{ \alpha
\in j\al_i+L\,|\, \la \alpha,\alpha\ra = 2\}$ be the set of all
roots in the coset $j\al_i+L$ for $j=1, \ldots, n_i-1$. Then
\[
\Phi(E_8)=\Phi \cup \bigcup_{j=1}^{n_i-1} H_j.
\]

We introduce weight $2$ elements $X^j$, namely,
\begin{equation}\label{eq:HWV}
X^j = \sum_{\alpha \in H_j} e^{\sqrt{2}\alpha}, \quad j=1, \ldots,
n_i-1.
\end{equation}
Then
\begin{equation}\label{eq:EF2}
\begin{split}
\hat{e} &= \frac{1}{16}\om +\frac{1}{32}\Big( \sum_{\al\in\Phi}
e^{\sqrt{2}\al}+ \sum_{j=1}^{n_i-1} X^j\Big),\\
\hat{f} &= \frac{1}{16}\om +\frac{1}{32}\Big( \sum_{\al\in\Phi}
e^{\sqrt{2}\al}+ \sum_{j=1}^{n_i-1} \xi^j X^j\Big),
\end{split}
\end{equation}
where $\xi = e^{2\pi\sqrt{-1}/n_i}$ is a primitive $n_i$-th root
of unity.

\begin{lemma}
$(1)$\ $X^j \in U$, $j = 1, \ldots , n_i-1$.

$(2)$\ $\hat{e}, \hat{f} \in U$.
\end{lemma}

\begin{proof}
Let $s^k$ be defined as in \eqref{eq:CV}. Then by a similar
argument as in the proof of Proposition \ref{hw}, we can verify
that $(s^k)_1 X^j = 0$ and $(s^k)_1 \hat{e} = 0$ for $k=1,\dots,
l$. Thus $X^j, \hat{e} \in U$ by the definition \eqref{eq:UDEF} of
$U$. Since $\sigma$ leaves $U$ invariant, we also have $\hat{f}
\in U$.
\end{proof}

\begin{remark}
The Weyl group $W(E_8)$ of the root system of type $E_8$ acts
naturally on the lattice VOA $V_{\sqrt{2}E_8}$ and $\hat{e}$ is
the only conformal vector among $\hat{e}^0_\delta,
\hat{e}^1_\zeta$ which is fixed by $W(E_8)$. The conformal vector
$\hat{f}$ is fixed by the Weyl group $W(\Phi) = W(\Phi_1) \times
\cdots \times W(\Phi_l)$ of the root system $\Phi = \Phi_1 \oplus
\cdots \oplus \Phi_l$ of $L=L(i)$. The conformal vector $\hat{e}$
is also fixed by the automorphism $\theta$ (cf. \eqref{eq:theta}).
However, $\hat{f}$ is not fixed by $\theta$ in general.
\end{remark}

\begin{theorem}
Let $\hat{e},\hat{f}$ be defined as in \eqref{eq:EF}. Then
\begin{equation}
\la \hat{e}, \hat{f}\ra = {\displaystyle
\begin{cases}
1/4 & \text{ if } i=0,\\
1/32 & \text{ if } i=1,\\
13/2^{10} & \text{ if } i=2,\\
1/2^{7} & \text{ if } i=3,\\
3/2^{9} & \text{ if } i=4,\\
5/2^{10} & \text{ if } i=5,\\
1/2^{8} & \text{ if } i=6,\\
0               & \text{ if } i=7,\\
1/2^{8} & \text{ if } i=8.
\end{cases}}
\end{equation}
In other words, the values of $\la \hat{e}, \hat{f}\ra$ are
exactly the values given in McKay's diagram \eqref{mckay}.
\end{theorem}

\begin{proof}
By \eqref{eq:EF2}, we can easily obtain that
\begin{equation*}
\la \hat{e}, \hat{f}\ra = \frac{1}{2^6}+ \frac{1}{2^{10}} \Big(
|\Phi| + \sum_{j=1}^{n_i-1} \xi^j |H_j|\Big),
\end{equation*}
where $H_j = \{ \alpha \in j\al_i+L\,|\, \la \alpha,\alpha\ra =
2\}$.

If $i=0$, then $n_0=1$ and $|\Phi|=240$. Hence
\begin{equation*}
\la\hat{e},
\hat{f}\ra=\frac{1}{2^6}+\frac{240}{2^{10}}=\frac{1}{4}.
\end{equation*}

If $i=1$, then $n_1=2$, $|\Phi|=|\Phi (A_1)|+|\Phi (E_7)|=128$,
and $| H_1|=112$. Hence
\begin{equation*}
\la\hat{e},
\hat{f}\ra=\frac{1}{2^6}+\frac{1}{2^{10}}(128-112)=\frac{1}{32}.
\end{equation*}

If $i=2$, then $n_2=3$, $|\Phi|=|\Phi (A_2)|+|\Phi (E_6)|=78$, and
$| H_1|=| H_2|=81$. Hence
\begin{equation*}
\la\hat{e},
\hat{f}\ra=\frac{1}{2^6}+\frac{1}{2^{10}}(78-81)=\frac{13}{2^{10}}.
\end{equation*}

If $i=3$, then $n_3=4$, $|\Phi|=|\Phi (A_3)|+|\Phi (D_5)|=52$,
$|H_1|=|H_3|=64$, and $| H_2|=60$. Hence
\begin{equation*}
\la\hat{e},
\hat{f}\ra=\frac{1}{2^6}+\frac{1}{2^{10}}(52-60)=\frac{1}{2^{7}}.
\end{equation*}

If $i=4$, then $n_4=5$, $|\Phi|=|\Phi (A_4)|+|\Phi (A_4)|=40$, and
$| H_1|=| H_2|=| H_3|=| H_4|=50$. Hence
\begin{equation*}
\la\hat{e},
\hat{f}\ra=\frac{1}{2^6}+\frac{1}{2^{10}}(40-50)=\frac{3}{2^{9}}.
\end{equation*}

If $i=5$, then $n_5=6$, $|\Phi|=|\Phi (A_1)|+|\Phi (A_2)|+|\Phi
(A_5)|=38$, $| H_1|=| H_5|=36$, $| H_2|=| H_4|=45$, and $|
H_3|=40$. Hence
\begin{equation*}
\la\hat{e},
\hat{f}\ra=\frac{1}{2^6}+\frac{1}{2^{10}}(38+36-45-40)=\frac{5}{2^{10}}.
\end{equation*}

If $i=6$, then $n_6=4$, $|\Phi|=|\Phi (A_1)|+|\Phi (A_7)|=58$, $|
H_1|=| H_3|=56$, and $| H_2|=70$. Hence
\begin{equation*}
\la\hat{e},
\hat{f}\ra=\frac{1}{2^6}+\frac{1}{2^{10}}(58-70)=\frac{1}{2^{8}}.
\end{equation*}

If $i=7$, then $n_7=2$, $|\Phi|=|\Phi (D_8)|=112$, and
$|H_1|=128$. Hence
\begin{equation*}
\la\hat{e}, \hat{f}\ra=\frac{1}{2^6}+\frac{1}{2^{10}}(112-128)=0.
\end{equation*}

If $i=8$, then $n_8=3$, $|\Phi|=|\Phi (A_8)|=72$, and
$|H_1|=|H_2|=84$. Hence
\begin{equation*}
\la\hat{e},
\hat{f}\ra=\frac{1}{2^6}+\frac{1}{2^{10}}(72-84)=\frac{1}{2^{8}}.
\end{equation*}

Thus we have proved the theorem.
\end{proof}

\begin{remark}
The same result still holds if we replace $\hat{e}$ by
$\hat{e}^\ep_\delta$ and $\hat{f}=\sigma{\hat{e}}$ by
$\sigma{\hat{e}^\ep_\delta}$ for any $\ep=0,1$ and  $\delta\in
{\Z_2}^8$.
\end{remark}

\section{Miyamoto's $\tau$-involutions and the canonical
automorphism $\sigma$}

Let $V$ be a VOA. If $V$ contains a conformal vector $w$ of
central charge $1/2$ such that the subalgebra $\Vir(w)$ generated
by $w$ is isomorphic to the Virasoro VOA $L(1/2,0)$, then an
automorphism $\tau_w$ of $V$ with $(\tau_w)^2 = 1$ can be defined.
Indeed, $V$ is a direct sum of irreducible $\Vir(w)$-modules.
Denote by $W_h$ the sum of all irreducible direct summands which
are isomorphic to $L(1/2,h)$. Then $\tau_w$ is defined to be $1$
on $W_0 \oplus W_{1/2}$ and $-1$ on $W_{1/16}$ (cf.
\cite{m1,m4}). Thus $\tau_w$ is the identity if $V$ has no
irreducible direct summand isomorphic to $L(1/2,1/16)$. We call
$\tau_w$ the Miyamoto involution or the $\tau$-involution
associated with $w$.

In this section, we shall study the relationship between the
canonical automorphism $\sigma$ and the Miyamoto involutions
$\tau_{\hat{e}}$, $\tau_{\sigma\hat{e}},\ldots$, and
$\tau_{\sigma^{n_i-1}\hat{e}}$. Let us recall two conformal
vectors $\hat{e}$ and $\hat{f}$ of central charge $1/2$ defined by
\eqref{eq:EF} and two automorphisms $\sigma$ and $\theta$
introduced in Subsection 3.1.

\begin{lemma}\label{tau-E8}
As automorphisms of $V_{\sqrt{2}E_8}$, $\tau_{\hat{e}}=\theta$.
\end{lemma}

\begin{proof}
By Proposition \ref{frame}, we know that $\{\hat{e}^0_\delta,
\hat{e}^1_\zeta\, |\, \delta, \zeta \in {\Z_2}^8/H_8,\ \delta,
\zeta \text{ even}\}$ is a Virasoro frame of $V_{\sqrt{2}E_8}^+$
and with respect to this frame, $V_{\sqrt{2}E_8}^+$ is a code VOA
isomorphic to $M_{RM(4,2)}$. Therefore, $\tau_{\hat{e}}|_
{V_{\sqrt{2}E_8}^+}= \mathrm{id}$. On the other hand,
\begin{equation*}
\begin{split}
\hat{e}_1 \gamma(-1)\cdot 1 &= \frac{1}{16}\omega_1
\gamma(-1)\cdot 1+\frac{1}{32}\sum_{\alpha\in\Phi^{+}(E_8)}
(e^{\sqrt{2}\alpha}+e^{-\sqrt{2}\alpha})_1 \gamma(-1)\cdot 1\\
& = \frac{1}{16} \gamma(-1)\cdot 1
\end{split}
\end{equation*}
for any $\gamma\in \sqrt{2}E_8$. By the definition of
$\tau_{\hat{e}}$, this implies that $\tau_{\hat{e}}
(\gamma(-1)\cdot 1 )=-\gamma(-1)\cdot 1$. Then $\tau_{\hat{e}}|_
{V_{\sqrt{2}E_8}^-}= -\mathrm{id}$, since $V_{\sqrt{2}E_8}^-$ is
an irreducible $V_{\sqrt{2}E_8}^+$-module. Hence
$\tau_{\hat{e}}=\theta$ as automorphisms of $V_{\sqrt{2}E_8}$.
\end{proof}

\begin{theorem}\label{tau-ef}
As automorphisms of $V_{\sqrt{2}E_8}$,
$\tau_{\hat{e}}\tau_{\hat{f}} = (\sigma^{-1})^2 = e^{2\pi
\sqrt{-1} \bbe(0)}$ and thus $|\tau_{\hat{e}}\tau_{\hat{f}}| =
n_i$ if $n_i$ is odd and $|\tau_{\hat{e}}\tau_{\hat{f}}| = n_i/2 $
if $n_i$ is even.
\end{theorem}

\begin{proof}
Since $\hat{f}=\sigma\hat{e}$, we have
$\tau_{\hat{f}}=\sigma\tau_{\hat{e}} \sigma^{-1}$. By \eqref{au}
and the preceding lemma, we also have
$\tau_{\hat{e}}\sigma\tau_{\hat{e}} =
\theta\sigma\theta=\sigma^{-1}$. Hence the assertion holds by
\eqref{sigmadef}.
\end{proof}

Next, we shall extend $\tau_{\hat{e}}, \tau_{\hat{f}}$, and
$\sigma$ to the Leech lattice VOA $V_{\Lambda}$. According to the
presentation \eqref{eq:2E8} of $\sqrt{2}E_8$, the dual lattice
$\mathcal{L}$ of $\sqrt{2}E_8$ is given by
\begin{equation*}
\mathcal{L} = \{ (a_1, \dots, a_8)\in \frac{1}2 \Z^8\, |\ 2(a_1,
\dots, a_8)\in H_8 \mod 2\}.
\end{equation*}

Note that $|\mathcal{L}/ \sqrt{2}E_8|=2^{8}$. Note also that
\[
V_{\mathcal{L}}= S(\mathfrak{h}^-_\Z)\otimes \C\{\mathcal{L}\}
\cong \bigoplus_{\al+\sqrt{2}E_8\in \mathcal{L}/ \sqrt{2}E_8}
V_{\al+\sqrt{2}E_8}
\]
as a module of $V_{\sqrt{2}E_8}$.

For any coset $\al+\sqrt{2}E_8$ of $\sqrt{2}E_8$ in $\mathcal{L}$,
one can always find a coset representative $\al$ whose square norm
is minimum in the coset such that $\al$ is in one of the following
forms.
\begin{equation}\label{eq:MCOSET}
\begin{split}
 &(0^8),\quad  (1, 0^7), \quad(1^2, 0^6), \quad  ((1/2)^4, 0^4), \\
 &((1/2)^3, -1/2, 0^4), \quad
   ((1/2)^2, (-1/2)^2, 0^4), \quad
   ((1/2)^4, 1, 0^3),\\
&((1/2)^3, -1/2,1, 0^3), \quad ((1/2)^8), \quad ((1/2)^7,
-1/2),\quad ((1/2)^6, (-1/2)^2).
\end{split}
\end{equation}

The square norm $\la \al,\al\ra$ of such $\al$ is $0$, $1$, or
$2$. Moreover, if $\la \al,\al\ra =2$, then $\al$ can be written
as a sum $ \al=a+b$, where $a,b \in \mathcal{L}$ are in the above
forms with $\la a,a \ra=\la b,b\ra=1$ and $\la a,b\ra=0$. In
particular, the minimal weight of the irreducible module
$V_{\al+\sqrt{2}E_8}$ is either $1/2$ or $1$ for $\alpha \not\in
\sqrt{2}E_8$.

Now $\sigma= e^{-\pi \sqrt{-1} \bbe(0)}$ (cf. \eqref{sigmadef})
acts on $V_{\mathcal{L}}$ as an automorphism of order $2n_i$. The
$\tau$-involution $\tau_{\hat{e}}$ also acts on $V_\mathcal{L}$.
In fact, $V_{\al+\sqrt{2}E_8}$ is $\tau_{\hat{e}}$-invariant for
any coset $\al +\sqrt{2}E_8$ of $\sqrt{2}E_8$ in $\mathcal{L}$.

\begin{lemma}\label{tau-ex}
For any $x\in \mathcal{L}$ with $\la x,x \ra =1$, $\tau_{\hat{e}}
(e^x) = -  e^{-x}$.
\end{lemma}

\begin{proof}
If $\la \gamma, \gamma \ra =4$ and $\la \gamma+x,\gamma+x \ra=1$
for some $\gamma \in \sqrt{2}E_8$, then $\la \gamma, x \ra=-2$ and
$\gamma+x=-x$. Thus, by the definition of $\hat{e}$ it follows
that
\[
\hat{e}_1 e^x = \frac{1}{16}\Big(\frac{1}2 e^x\Big) + \frac{1}{32}
e^{-x} \quad \text{and} \quad \hat{e}_1 e^{-x} =
\frac{1}{16}\Big(\frac{1}2 e^{-x} \Big) + \frac{1}{32} e^{x}.
\]
Therefore, $\hat{e}_1(e^x+e^{-x})= \frac{1}{16} (e^x+e^{-x})$ and
$\hat{e}_1(e^x- e^{-x}) =0$. Hence $\tau_{\hat{e}} (e^x+e^{-x}) =
- (e^x+e^{-x})$ and $\tau_{\hat{e}} (e^x-e^{-x}) =e^x-e^{-x}$ by
the definition of $\tau_{\hat{e}}$, and so $\tau_{\hat{e}} (e^x) =
-  e^{-x}$.
\end{proof}

\begin{lemma}
Let $\al+\sqrt{2}E_8$ be a coset of $\sqrt{2}E_8$ in
$\mathcal{L}$. Then for any $u\in V_{\al+\sqrt{2}E_8}$,
$\tau_{\hat{e}} \sigma \tau_{\hat{e}} (u) = \sigma^{-1}(u)$.
\end{lemma}

\begin{proof}
We have $V_{\al+\sqrt{2}E_8}= \mathrm{span}_{\C}\{ v_n e^\al\,|\,
v\in V_{\sqrt{2}E_8},\,n\in \Z\}$, since $V_{\al+\sqrt{2}E_8}$ is
an irreducible $V_{\sqrt{2}E_8}$-module. If $\la \al,\al\ra=1 $,
then we know that $\tau_{\hat{e}}( e^\al)=-e^{-\al}$ by Lemma
\ref{tau-ex}. Thus $\tau_{\hat{e}} \sigma \tau_{\hat{e}} (e^\al) =
\sigma^{-1}(e^\al)$ and so
\[
\begin{split}
\tau_{\hat{e}} \sigma \tau_{\hat{e}} (v_n e^\al )
&=(\tau_{\hat{e}} \sigma \tau_{\hat{e}} (v))_n\,  (\tau_{\hat{e}}
\sigma \tau_{\hat{e}}(e^\al))\\
&= \sigma^{-1}(v)_n \sigma^{-1} (e^\al)\\
&= \sigma^{-1} (v_n e^\al)
\end{split}
\]
for any $v \in V_{\sqrt{2}E_8}$ by Lemma \ref{tau-E8}.

If $\la \al,\al\ra=2 $, then $\al=a+b$ for some vectors $a, b$ in
the forms of \eqref{eq:MCOSET} with $\la a,a \ra=\la b,b\ra=1$ and
$\la a,b\ra=0$. In this case, $e^\al= (e^a)_{-1} e^b$ and we still
have $\tau_{\hat{e}} \sigma \tau_{\hat{e}} (e^\al) =
\sigma^{-1}(e^\al)$. Thus for any $v\in V_{\sqrt{2}E_8}$,
\[
\begin{split}
\tau_{\hat{e}} \sigma \tau_{\hat{e}} (v_n e^\al ) =
\sigma^{-1}(v)_n \sigma^{-1} (e^\al) = \sigma^{-1} (v_ne^\al)
\end{split}
\]
as required.
\end{proof}

As a consequence,  we have the following proposition.

\begin{proposition}\label{dual}
For any  $u\in V_{\mathcal{L}}$, $\tau_{\hat{e}} \sigma
\tau_{\hat{e}} (u) = \sigma^{-1}(u)$. Hence $\tau_{\hat{e}}
\tau_{\hat{f}}= (\sigma^{-1})^2= e^{2\pi \sqrt{-1} \bbe(0)}$ as
automorphisms of $V_{\mathcal{L}}$.
\end{proposition}

Now we discuss the situation in the Leech lattice VOA $V_\Lambda$.
First let us recall the following theorem \cite[Theorem 4.1]{dlmn}
(see also \cite{hk,lm}).
\begin{theorem}
For any even unimodular lattice $N$ of rank $24$, there is at
least one (in general many) isometric embedding of $\sqrt{2}N$
into the Leech lattice $\Lambda$.
\end{theorem}

It is well known (cf. \cite{hk}) that the Leech lattice $\Lambda$
can be constructed by ``Construction A" for $\Z_4$-codes of length
$24$. In fact,
\[
\Lambda = A_4(\mathcal{C})=\frac{ 1} 2 \{x\in \Z^{24}\,|\,x \equiv
c \mod 4\quad \text{ for some } c\in \mathcal{C}\}
\]
for some type II self-dual $\Z_4$-code $\mathcal{C}$ of length
$24$. By \cite{hk}, $\mathcal{C}$ can be taken to be the
$\Z_4$-code having the generating matrix \eqref{eq:GM4}.

\begin{equation}\label{eq:GM4}
\begin{pmatrix}
2222 &0000 &0000 &0000 &0000 &0000\\
0022 &2200 &0000 &0000 &0000 &0000 \\
0000 &0022 &2020 &0000 &0000 &0000 \\
0000 &0000 &0202 &2020 &0000 &0000 \\
0000 &0000 &0000 &0202 &2002 &0000 \\
2020 &2020 &0000 &0000 &0000 &0000 \\
0000 &0220 &2200 &0000 &0000 &0000 \\
0000 &0000 &2002 &2002 &0000 &0000 \\
0000 &0000 &0000 &0022 &2020 &0000 \\
2000 &2000 &2000 &2000 &2000 &2000 \\
1111 &1111 &2000 &2000 &0000 &0000 \\
2000 &1111 &1111 &0000 &2000 &0000 \\
0000 &0000 &1111 &1111 &2000 &2000 \\
2000 &0000 &2000 &1111 &1111 &0000 \\
2000 &2000 &0000 &0000 &1111 &1111 \\
3012 &1010 &1001 &1001 &1100 &1100 \\
3201 &1001 &1100 &1100 &1010 &1010
\end{pmatrix}
\end{equation}

For any $\Z_4$-code $C$ of length $n$, one can obtain a binary
code
\[
B(C)=\{ (b_1, \dots,b_n) \in {\Z_2}^n\,|\, (2b_1,\dots,2b_n)\in
C\},
\]
where $2b_j$ should be considered as $0 \in \Z_4$ if $b_j = 0 \in
\Z_2$ and $2 \in \Z_4$ if $b_j = 1 \in \Z_2$.  Moreover, the
lattice
\[
L_{B(C)}=\{ x\in \Z^n\,|\ x\in B(C) \mod 2\}
\]
is a sublattice of $A_4(C)$. In the case for $C=\mathcal{C}$, the
binary code $B(\mathcal{C})$ contains a subcode isomorphic to
$H_8\oplus H_8\oplus H_8$. Thus by \eqref{eq:2E8}, we have an
explicit embedding of $\sqrt{2}{E_8}^3$ into the Leech lattice
$\Lambda$.

Now let $ \sqrt{2}{E_8}^3 \longrightarrow \Lambda$ be any
embedding of $\sqrt{2}{E_8}^3$ into the Leech lattice $\Lambda
\subset \mathcal{L}^3$. Let $\tilde{\bbe}=\sqrt{2}(\ba,0,0) \in
\mathcal{L}^3$, where $\ba$ is defined as in \eqref{adef}. Define
$\tilde{\sigma} : (V_\mathcal{L})^{\otimes 3} \rightarrow
(V_\mathcal{L})^{\otimes 3}$ by
\[
\tilde{\sigma}=\sigma\otimes 1\otimes 1 = e^{-\pi \sqrt{-1}
\tilde{\bbe}(0)}.
\]
Then $\tilde{\sigma}$ is an automorphism of $V_{\Lambda}$.
Moreover, the following theorem holds.

\begin{theorem}
Let $\tilde{\bbe}$ and $\tilde{\sigma}$ be defined as above. Then
as automorphisms of $V_{\Lambda}$, $\tau_{\hat{e}}\tau_{\hat{f}} =
(\tilde{\sigma}^{-1})^2= e^{2\pi \sqrt{-1} \tilde{\bbe}(0)}$ and
$|\tau_{\hat{e}}\tau_{\hat{f}}|=n_i$  for any $i=0,1,\dots,8$.
\end{theorem}

\section{Correspondence with Conway's axes.}
Recall the elements $\tilde{\omega}^k$ and $X^j$ defined by
\eqref{eq:CV} and \eqref{eq:HWV}. It turns out that the Griess
algebra $U_2$ of $U$ is generated by $\hat{e}$ and $\hat{f}$ and
is of dimension $l + n_i -1$ with basis $\tilde{\omega}^k, 1 \le k
\le l$ and $X^j, 1 \le j \le n_i-1$ (see \cite{lyy} for details).
We can verify that the Griess algebra $U_2$ coincides with the
algebra described in Conway \cite[Table 3]{c1}. In \cite{c1}, it
is shown that for each $2A$-involution of the Monster simple
group, there is a unique idempotent in the Monstrous Griess
algebra $V^\natural_2$ corresponding to the involution. Such an
idempotent is called an axis. By Miyamoto \cite{m1}, an axis is
exactly half of a conformal vector of central charge $1/2$. Note
that the product $t \ast t'$ and the inner product $(t,t')$ of two
axes $t,t'$ in \cite{c1} are equal to $t\cdot t'=t_1 t'$ and $\la
t,t'\ra/2$, respectively in our notation. Let $t_n$ be as in
\cite{c1}. We denote $t$, $u$, $v$, and $w$ of \cite{c1} by
$t_{2A}$, $u_{3A}$, $v_{4A}$, and $w_{5A}$, respectively.

In each of the nine cases, we obtain an isomorphism of our Griess
algebra $U_2$ to Conway's algebra generated by two axes through
the following correspondence between our conformal vectors and
Conway's axes.

\medskip
$1A$ case.\qquad  $\hat{e} \longleftrightarrow \frac{1}{32}t_0$.

\medskip
$2A$ case.\qquad $\sigma^j \hat{e} \longleftrightarrow
\frac{1}{32}t_j,\ j=0,1$,\quad $\tilde{\omega}^1
\longleftrightarrow \frac{1}{32}t_{2A}$.

\medskip
$3A$ case.\qquad $\sigma^j \hat{e} \longleftrightarrow
\frac{1}{32}t_j,\ j=0,1,2$,\quad $\tilde{\omega}^1
\longleftrightarrow \frac{1}{45}u_{3A}$.

\medskip
$4A$ case.\qquad $\sigma^j \hat{e} \longleftrightarrow
\frac{1}{32}t_j,\ 0 \le j \le 3$,\quad $\tilde{\omega}^1
\longleftrightarrow \frac{1}{96}v_{4A}$.

\medskip
$5A$ case.\qquad $\sigma^j \hat{e} \longleftrightarrow
\frac{1}{32}t_j,\ 0 \le j \le 4$,\quad $\tilde{\omega}^1 -
\tilde{\omega}^2 \longleftrightarrow -\frac{1}{35\sqrt{5}}w_{5A}$.

\medskip
$6A$ case.\qquad $\sigma^j \hat{e} \longleftrightarrow
\frac{1}{32}t_j,\ 0 \le j \le 5$,\quad $\tilde{\omega}^2
\longleftrightarrow \frac{1}{32}t_{2A}$, \quad $\tilde{\omega}^1
\longleftrightarrow \frac{1}{45}u_{3A}$.

\medskip
$4B$ case.\qquad $\sigma^j \hat{e} \longleftrightarrow
\frac{1}{32}t_j,\ 0 \le j \le 3$,\quad $\tilde{\omega}^1
\longleftrightarrow \frac{1}{32}t_{2A}$.

\medskip
$2B$ case.\qquad $\sigma^j \hat{e} \longleftrightarrow
\frac{1}{32}t_j,\ j=0,1$.

\medskip
$3C$ case.\qquad $\sigma^j \hat{e} \longleftrightarrow
\frac{1}{32}t_j,\ j=0,1,2$.

\bibliographystyle{amsplain}

\begin{thebibliography}{99}
\bibitem{c1}
J. H. Conway, A simple construction for the Fisher-Griess
Monster group, \emph{Invent. Math}. \textbf{79} (1985), 513--540.

\bibitem{atlas}
J. H. Conway, R. T. Curtis, S. P. Norton, R. A. Parker and R. A.
Wilson, \emph{ATLAS of Finite Groups}, Oxford Univ. Press, 1985.

\bibitem{cs}
J. H. Conway and N. J. A. Sloane, \emph{Sphere Packings, Lattices
and Groups}, Springer-Verlag, Berlin-New York, 1988.

\bibitem{d1}
C. Dong, Vertex algebras associated with even lattices, \emph{J.
Algebra} \textbf{161} (1993), 245--265.

\bibitem{dlmn}
C. Dong, H. Li, G. Mason and S. P. Norton, Associative subalgebras
of Griess algebra and related topics, \emph{Proc. of the
Conference on the Monster and Lie algebra at the Ohio State
University}, May 1996, ed. by J. Ferrar and K. Harada, Walter de
Gruyter, Berlin-New York, 1998, pp. 27--42.

\bibitem{fhl}
I. B. Frenkel, Y. Huang and J. Lepowsky, \emph{On axiomatic
approaches to vertex operator algebras and modules}, Mem. Amer.
Math. Soc. \textbf{104}, 1993.

\bibitem{flm}
I. B. Frenkel, J. Lepowsky and A. Meurman,\emph{Vertex Operator
Algebras and the Monster}, Pure and Applied Math., Vol.
\textbf{134}, Academic Press, New York, 1988.

\bibitem{GN}
G. Glauberman and S. P. Norton, On McKay's connection between
the affine $E_8$ diagram and the Monster, \emph{CRM Proceedings
and Lecture Notes}, Vol. \textbf{30}, Amer. Math. Soc.,
Providence, RI, 2001, pp. 37--42.

\bibitem{griess}
R. Griess, The Friendly Giant, \emph{Invent. Math.} \textbf{69}
(1982), 1--102.

\bibitem{hk}
M. Harada and M. Kitazume, $\Z_4$-code constructions for the
Niemeier lattices and their embeddings in the Leech lattice,
\emph{Europ. J. Combinatorics } \textbf{21} (2000), 473--485.

\bibitem{ly3}
C. H. Lam and H. Yamada, Decomposition of the lattice vertex
operator algebra $V_{\sqrt{2}A_l}$, \emph{J. Algebra} \textbf{272}
(2004), 614--624.

\bibitem{lyy}
C. H. Lam, H. Yamada and H. Yamauchi, McKay's observation and
vertex operator algebras generated by two conformal vectors of
central charge $1/2$, in preparation.

\bibitem{lm}
J. Lepowsky and A. Meurman, An $E_8$-approach to the Leech lattice
and the Conway group, \emph{J. Algebra} \textbf{77} (1982),
484--504.

\bibitem{McK}
J. McKay, Graphs, singularities, and finite groups,
\emph{Proc. Symp. Pure Math.}, Vol. \textbf{37}, Amer. Math. Soc.,
Providence, RI, 1980, pp. 183--186.

\bibitem{m1}
M. Miyamoto, Griess algebras and conformal vectors in vertex
operator algebras, \emph{J. Algebra} \textbf{179} (1996),
523--548.

\bibitem{m2}
M. Miyamoto, Binary codes and vertex operator (super)algebras,
\emph{J. Algebra} \textbf{181} (1996), 207--222.

\bibitem{m4}
M. Miyamoto, A new construction of the Moonshine vertex operator
algebras over the real number field, to appear in \emph{Ann. of
Math.}

\bibitem{zf1}
A. B. Zamolodchikov and V. A. Fateev, Nonlocal (parafermion)
currents in two dimensional conformal quantum field theory and
self-dual critical points in $\Z_N$-symmetric statistical systems,
\emph{Sov. Phys. JETP} \textbf{62} (1985), 215--225.

\end{thebibliography}

\end{document}